\documentclass{amsart}
\usepackage{amssymb}
\usepackage{amscd,enumerate,amsfonts,calc,amsmath,verbatim}
\usepackage[all]{xy}
\usepackage{enumerate}
\usepackage{tikz}
\usepackage{xcolor}
\usepackage[hidelinks]{hyperref}
\hypersetup{colorlinks,breaklinks,
             urlcolor=blue,
             linkcolor=blue}
\usepackage{graphicx}
\DeclareGraphicsExtensions{.pdf,.jpeg,.png}
\newtheorem{thm}{{\bf Theorem}}
\newtheorem{lemma}[thm]{{\bf Lemma}}
\newtheorem{rmk}[thm]{{\bf Remark}}
\newtheorem{prop}[thm]{{\bf Proposition}}
\newtheorem{cor}[thm]{{\bf Corollary}}

\newcommand{\F}{\mathbb F}
\newcommand{\N}{\mathbb N}

\def\1{{1\hskip-0.25em{\rm 1}}}

\def\CC{{\ \rlap{\raise 0.4ex 
\hbox{$\scriptscriptstyle |$}}\hskip -0.2em C}}

\def\dim{\text{dim}}

\def\Fq{{\mathbb F}_q}

\newcommand{\Ev}{\operatorname{Ev}}
\newcommand{\supp}{\mathrm{Supp}}

\newcommand{\wt}{\mathrm{wt}}

\begin{document}
\title[Purity of  Resolutions Associated to Reed-Muller Codes]{On the Purity of  Resolutions of Stanley-Reisner Rings Associated to Reed-Muller Codes}
	\author{Sudhir R. Ghorpade}
	\address{Department of Mathematics, 
		Indian Institute of Technology Bombay,\newline \indent
		Powai, Mumbai 400076, India}
\email{\href{mailto:srg@math.iitb.ac.in}{srg@math.iitb.ac.in}}
\thanks{Sudhir Ghorpade is partially supported by DST-RCN grant INT/NOR/RCN/ICT/P-03/2018 from the Dept. of Science \& Technology, Govt. of India, MATRICS grant MTR/2018/000369 from the Science \& Engg. Research Board, and IRCC award grant 12IRAWD009 from IIT Bombay.}
	\author{Rati Ludhani}
	\address{Department of Mathematics, 
		Indian Institute of Technology Bombay,\newline \indent
		Powai, Mumbai 400076, India}
\email{\href{mailto:lrati@math.iitb.ac.in}{lrati@math.iitb.ac.in}}
\thanks{Rati Ludhani is supported by Prime Minister's Research Fellowship PMRF-192002-256 at IIT Bombay.}

\begin{abstract}
Following Johnsen and Verdure (2013), we can associate to any linear code $C$ an abstract simplicial complex and in turn, a Stanley-Reisner ring $R_C$. The ring $R_C$ is a standard graded algebra over a field and its projective dimension is precisely the dimension of $C$. Thus $R_C$ admits a graded minimal free resolution and the resulting graded Betti numbers are known to determine the generalized Hamming weights of $C$.  The question of purity of the minimal free resolution of $R_C$ was considered by Ghorpade and Singh (2020) when $C$ is the generalized Reed-Muller code. They showed that the resolution is pure in some cases and it is not pure in many other cases. Here we give a complete characterization of the purity of graded minimal free resolutions of Stanley-Reisner rings associated to generalized Reed-Muller codes of an arbitrary order. 
\end{abstract}

\maketitle
\section{introduction}

This article concerns a topic that is at the interface of homological aspects of commutative algebra and the theory of linear error correcting codes. Our motivation comes from the work of Johnsen and Verdure \cite{JV1} and the more recent work \cite{GS}. In \cite{JV1}, 
 the notion of {\it Betti numbers} of a linear code is introduced. The Betti numbers of a linear code $C$ of length $n$ are, in fact, the graded Betti numbers of the Stanley-Reisner ring $R_C$ of the simplicial complex $\Delta_C$ on $[n]:=\{1, \dots , n\}$ whose faces are precisely the subsets $\{i_1, \dots , i_t\}$ of $[n]$ for which the columns $H_{i_1}, \dots , H_{i_t}$ of a parity check matrix $H$ of $C$ are linearly independent. 
 In  \cite{JV1}, it was shown that the Betti numbers of a linear code determine its generalized Hamming weights. 
 Further, \cite{JRV, JP} showed that the Betti numbers of a linear code (and its elongations) are also  closely related to several classical parameters of that code. Thus it is useful to know them explicitly.  Computation of these Betti numbers is in general, a difficult problem, but it becomes easy, by a formula of Herzog and K{\"{u}}hl  \cite{HK}, when the corresponding minimal free resolutions are pure. An intrinsic characterization of purity of the graded minimal free resolutions of Stanley-Reisner rings associated to arbitrary linear codes was obtained in \cite{GS}. As a consequence, known results about the Betti numbers 
 of  MDS codes and constant weight codes were easily deduced. 
 
 One of the most important and widely studied class of linear codes is that of Reed-Muller codes. These codes were introduced by Reed \cite{R} in the binary case and several of their properties were established by Muller \cite{M}; see also \cite[pp. 20--38]{B}. We shall consider Reed-Muller codes in the most general sense, as given by Kasami, Lin and Peterson \cite{KLP} and by Delsarte, Goethals, and MacWilliams \cite{DGW}. Generalized Hamming weights of (generalized) Reed-Muller codes are explicitly known, thanks to the work of Heijnen and Pellikaan \cite{HP} (see also \cite{BD} and \cite{Be}).  It is, therefore, natural, to ask for an explicit determination of the Betti numbers of Reed-Muller codes. The problem would be tractable if we know when the graded minimal free resolutions of Stanley-Reisner rings of simplicial complexes corresponding to Reed-Muller codes are pure. This question about purity was considered in \cite{GS} and an answer was provided in many, but not all, cases. In this article we build upon the work in \cite{GS} and complete it  to give 
  a 
  characterization of purity of  graded minimal free resolutions of Stanley-Reisner rings associated to arbitrary Reed-Muller codes. 
 
 This paper is organized as follows. In Section~\ref{sec2}, 
we review (generalized) Reed-Muller codes  and discuss their properties that are relevant for us. 
 Next, in Section~\ref{sec3}, the notion of purity of a minimal free resolution is recalled and some key results in \cite{GS} such as the intrinsic characterization mentioned above and 
 results about the purity or non-purity of resolutions corresponding to Reed-Muller codes are stated. Our main result on a 
 characterization of  purity of 
 free resolutions of Stanley-Reisner rings associated to Reed-Muller codes is also proved here. 
As a corollary, we  give a characterization of Reed-Muller codes that are MDS codes.

\section{Reed-Muller codes}\label{sec2}

Standard references for (generalized) Reed-Muller codes are the book of Assmus and Key \cite{AK} (especially, Chapter~5) and the seminal paper of Delsarte, Goethals, and MacWilliams \cite{DGW}. 
Let us begin by setting some basic notation and terminology. 

Fix throughout this paper a prime power $q$ and a finite field $\F_q$ with $q$ elements. Let $n,k$ be integers with $1\le k \le n$. We write \emph{$[n,k]_q$-code} to mean a $q$-ary linear code of length $n$ and dimension $k$, i.e., a $k$-dimensional $\Fq$-linear  subspace of $\Fq^n$. If the minimum distance of an $[n,k]_q$-code is $d$, then it may be referred to as an $[n,k,d]_q$-code. If $C$ is an $[n,k,d]_q$-code, then the elements of $C$ of Hamming weight $d$ will be referred to as the \emph{minimum weight codewords} of $C$. 
An $[n,k]_q$-code is said to be \emph{nondegenerate} if it is not contained in a coordinate hyperplane of $\Fq^n$. 
We denote by $\N$ the set of nonnegative integers. 

Let $m, r$ be integers such that $m\geq 1$ and $0\leq r\leq m(q-1)$. Define
\begin{eqnarray*}
V_q(r,m):=\{f\in \F_q[X_1,\ldots,X_m]:{\rm deg}(f)\leq r\;{\rm and}\; {\rm deg}_{X_i}(f)<q\; {\rm for}\; i=1,\ldots,m\}.
\end{eqnarray*}
Note that 
$V_q(r,m)$ is a $\Fq$-linear subspace of the polynomial ring $\F_q[X_1,\ldots,X_m]$. Fix an ordering 
${\mathsf{P}}_1,\ldots,{\mathsf{P}}_{q^m}$ of the elements of $\Fq^m$ and consider the evaluation map 
\begin{eqnarray}\label{eq:Ev}
\Ev:V_q(r,m) \to \Fq^{q^m} \quad {\rm defined\; by}\quad  f\mapsto c^{ }_f:=(f({\mathsf{P}}_1),\ldots, f({\mathsf{P}}_{q^m})).
\end{eqnarray}
Clearly, $\Ev$ is a linear map and its image is a nondegenerate linear code of length $q^m$; this code 
is called the \emph{(generalized) Reed-Muller code of order $r$}, and it is denoted by $\mathcal{RM}_q(r,m)$. 
The dimension of $\mathcal{RM}_q(r,m)$ is given by the following 
formula that can be found in Assmus and Key \cite[Theorem 5.4.1]{AK}: 
\begin{equation}\label{AKdim}
\dim\, \mathcal{RM}_q(r,m) =  \sum_{s=0}^r \sum_{i=0}^m (-1)^i \binom{m}{i} \binom{s-iq + m - 1}{s-iq}.
\end{equation}
In \cite[eq. (13)]{GS}, a somewhat simpler formula for the dimension is stated (without proof). 
It is not difficult to derive it from \eqref{AKdim}. However, we give an independent and direct proof of the simpler formula below. 

\begin{lemma}\label{lem1}
Let $m, r$ be integers such that $m\geq 1$ and $0\leq r\leq m(q-1)$. Then 
\begin{equation}\label{GSdim}
\mathop{\rm dim} \mathcal{RM}_q(r,m) = \sum_{i=0}^{m}(-1)^i \binom{m}{i}\binom{m+r-iq}{m}.
\end{equation}
\end{lemma}

\begin{proof}
It is well-known that the map $\Ev$ given by \eqref{eq:Ev} is injective. This follows, for instance, from \cite[Lemma 2.1]{G}. Also, if 
$E: = \{(v_1, \dots , v_m) \in \N^m : v_1+ \dots + v_m \le r \}$, then it is easily seen that a basis of $V_q(r,m)$ is given by 
$$
B:= \{X_1^{v_1} \cdots X_m^{v_m} : (v_1, \dots , v_m) \in  E \text{ and } 0 \le v_j < q \text{ for } 1\le j \le m\}.
$$
Let $E_j :=  \{(v_1, \dots , v_m) \in E : v_j \ge q\}$  for $1\le j \le m$. 
The set $B$ is clearly in bijection with $E \setminus (E_1 \cup \dots \cup E_m)$. 
It is elementary and well-known that $|E| = \binom{m+r}{m}$. By changing $v_j$ to $v'_j = v_j -q$, we also see that 
$|E_j| = {{m+r-q}\choose{m}}$ for $1\le j \le m$, and more generally, 
$|E_{j_1} \cap \dots \cap E_{j_i}| = {{m+r-iq}\choose{m}}$ for $1 \le j_1 < \dots < j_i \le m$. 
It follows that $\dim\, \mathcal{RM}_q(r,m) = \dim\, V_q(r,m) = |B|$,  and this is equal to 
\begin{eqnarray*}
|E| -|E_1 \cup \dots \cup E_m|
& = & \binom{m+r}{m} - \sum_{i=1}^{m} (-1)^{i-1} \sum_{1 \le j_1 < \dots < j_i \le m} |E_{j_1} \cap \dots \cap E_{j_i}| \\
& = & 
\binom{m+r}{m} -  \sum_{i=1}^{m}(-1)^{i-1} \binom{m}{i}\binom{m+r-iq}{m}.
\end{eqnarray*}
The last expression is clearly equal to the desired formula in \eqref{GSdim}. 
\end{proof}

\begin{rmk}\label{rem:RMdim}
{\rm 
In case $0\le r < q$, 
formula \eqref{GSdim} simplifies to $ \mathop{\rm dim} \mathcal{RM}_q(r,m) = \binom{m+r}{m}$. This can also be seen by noting that the set $E_j$ in the proof above is empty for each $j=1, \dots , m$ when $r< q$.  On the other hand, if 
$r = m(q-1)$, then the map $\Ev$ given by \eqref{eq:Ev} is 
also surjective. To see this, write ${\mathsf{P}}_{\nu} = (a_{\nu 1}, \dots , a_{\nu m})$ and consider 
\begin{equation}\label{eq:Fnu}
F_{\nu}(X_1, \dots , X_m) := \prod_{j=1}^m \left( 1 - (X_j - a_{\nu j})^{q-1}\right) \quad \text{for } \nu = 1, \dots , q^m.
\end{equation}
Note that for any $\nu \in \{1, \dots , q^m\}$, the polynomial $F_{\nu}$ is in $V_q(m(q-1), m)$ and it has the property that
$F_{\nu}({\mathsf{P}}_{\nu}) = 1$ and $F_{\nu}({\mathsf{P}}_{\mu}) =0$ for any $\mu \in \{1, \dots , q^m\}$ with $\mu \ne \nu$. Hence
any $\lambda = (\lambda_1, \dots , \lambda_{q^m}) \in \Fq^{q^m}$ can be written as 
$\lambda = \Ev(F)$, where 
$F = \lambda_{1} F_{1} + \cdots + \lambda_{q^m} F_{q^m}$. It follows that $\mathcal{RM}_q(m(q-1),m) = \Fq^{q^m}$. In particular, Lemma~\ref{lem1} yields the following curious identity:
$$
\sum_{i=0}^{m}(-1)^i \binom{m}{i}\binom{(m-i)q}{m} = q^m \quad \text{or equivalently,} \quad 
\sum_{i=0}^{m}(-1)^i \binom{m}{i}\binom{iq}{m} = (-q)^m.
$$
It may be interesting to obtain a direct proof of the above identity.
}
\end{rmk}

We now recall 
the following important result about the minimum distance and the minimum weight codewords of Reed-Muller codes. 

\begin{prop}\label{RM-min}
Let $m, r$ be integers such that $m\geq 1$ and $0\leq r\leq m(q-1)$. Then there are unique $t,s\in \N$ such that 
\begin{equation}\label{eq:ts}
r=t(q-1)+s \quad  \text{and} \quad 0\leq s\leq q-2.
\end{equation}
With $t, s$ as above, the minimum distance of  $\mathcal{RM}_q(r,m)$ is given by 
\begin{equation}\label{eq:dRM}
d=(q-s)q^{m-t-1}. 
\end{equation}
Further, if $f\in V_q(r,m)$ is given by 
\begin{eqnarray}\label{eq:f}
f(X_1,\ldots,X_m)=\omega_0\prod_{i=1}^t \left(1-(X_i-\omega_i)^{q-1}\right) \prod_{j=1}^s (X_{t+1}-\omega'_j)
\end{eqnarray}
where $\omega_0,\omega_1,\ldots,\omega_t\in \F_q$ with $\omega_0\neq 0$ and $\omega'_1,\ldots,\omega'_s$ are any distinct elements of $\F_q$, then $\Ev(f)$ is a minimum weight codeword of $\mathcal{RM}_q(r,m)$. 
 Moreover, every minimum weight codeword of $\mathcal{RM}_q(r,m)$ is of the form $\Ev(g)$, where $g$ is obtained from a polynomial of the form \eqref{eq:f} by substituting for $X_1, \dots , X_{t+1}$ any $(t+1)$ linearly independent linear forms in $\Fq[X_1, \dots , X_m]$. 
\end{prop}

\begin{proof}
The formula in \eqref{eq:dRM} follows from \cite[Theorem 2.6.1]{DGW} and \cite[Theorem 5]{KLP}. The assertion about the minimum weight codewords 
is proved in  \cite[Theorem 2.6.3]{DGW}. 
\end{proof}

We end this section by observing 
that the Reed-Muller code $ \mathcal{RM}_q(r, m)$ is a particularly nice code when $m$ is small or when $r$ is either very small or very large. 

\begin{lemma}\label{lem2}
Let $m, r$ be integers such that $m\geq 1$ and $0\leq r\leq m(q-1)$. Then $ \mathcal{RM}_q(r, m)$ is an MDS code in each of the following cases: (i) $m=1$, (ii) $r=0$, (iii) $r = m(q-1)$,  and (iv) $r= m(q-1)-1$.
\end{lemma}

\begin{proof}
(i) if $0\le r < q$, then in view of Remark~\ref{rem:RMdim} and Proposition~\ref{RM-min}, we see that $\mathcal{RM}_q(r, 1)$ is a $[q, \, r+1, \, q-r]_q$-code, and hence it is an MDS code. 

(ii) Clearly, $\mathcal{RM}_q(0, m)$ is the $1$-dimensional code of length $q^m$ spanned by the all-1 vector, and this is evidently  an MDS code. 

(iii) 
From Remark~\ref{rem:RMdim}, 
$\mathcal{RM}_q(m(q-1), m) = \Fq^{q^m}$, which is  obviously 
an MDS~code. 

(iv) Suppose $r = m(q-1)-1$. We will show that 
\begin{equation}\label{eq:oneless}
\mathcal{RM}_q(r, m) = \Lambda,  \ \text{ where } \ \Lambda:= \left\{(\lambda_1, \dots , \lambda_{q^m}) \in \Fq^{q^m} : \lambda_1 + \dots + \lambda_{q^m} =0\right\}. 
\end{equation}
This would imply that $\mathcal{RM}_q(r, m)$ is a $[q^m, q^m-1, 2]_q$-code, and hence an MDS code. To prove \eqref{eq:oneless}, 
first note that the monomial $X_1^{q-1} \cdots X_m^{q-1}$ is in $V_q(m(q-1), m)$, but not in the subspace $V_q(r, m)$. Since we have seen in Remark~\ref{rem:RMdim} that $\Ev$ gives an isomorphism of $V_q(m(q-1), m)$ onto $\Fq^{q^m}$, it follows that $\dim_{\Fq} V_q(r, m) \le q^m-1$. 
Hence  
it suffices to show that $\Lambda \subseteq \mathcal{RM}_q(r, m)$. 
To this end, we 
assume without loss of generality that the ordering ${\mathsf{P}}_1, \dots , {\mathsf{P}}_{q^m}$ of points of $\Fq^m$ is such that ${\mathsf{P}}_1$ is the origin. For $1\le \nu \le q^m$, consider the polynomial $F_{\nu}$ given by \eqref{eq:Fnu}, and write
$$
F_{\nu} = H + G_{\nu}, \quad \text{where} \quad H := F_1 = \prod_{j=1}^m \big( 1 - X_j^{q-1}\big) \quad \text{and}
\quad G_{\nu}:= F_{\nu}-H. 
$$
Note that $G_{\nu} \in V_q(r, m)$ for each  $ \nu = 1, \dots , q^m$. 
Also, $H({\mathsf{P}}_1)=1$ and $H({\mathsf{P}}_{\mu}) =0$ for $2\le \mu \le q^m$. So in view of the properties of $F_{\nu}$ noted in Remark~\ref{rem:RMdim}, we see that 
$G_1({\mathsf{P}}_1)=0$ while $G_{\nu}({\mathsf{P}}_1) = -1$ and $G_{\nu}({\mathsf{P}}_{\nu})=1$ for $2\le \nu \le q^m$, and moreover, $G_{\nu}({\mathsf{P}}_{\mu})=0$ for $2\le \nu, \mu  \le q^m$ with $\nu \ne \mu$. 
Thus 
given any $\lambda =(\lambda_1, \dots , \lambda_{q^m}) \in \Lambda$, 
the polynomial $G:=\sum_{\lambda=1}^{q^m}\lambda_{\nu}G_{\nu} \in V_q(r,m)$ and 
$
\Ev(G) = \lambda. 
$
This proves \eqref{eq:oneless}. 
\end{proof}

\begin{rmk}\label{rem:RMsmall}
{\rm
In \cite[pp. 8--9]{GS}, the results in Lemma~\ref{lem2}, especially (iv), were deduced by appealing to the structure of duals of Reed-Muller codes. Here we have chosen to give a more direct and elementary proof. We remark also that the converse of the result in Lemma~\ref{lem2} is true. An indirect proof of this is given later; see  Corollary~\ref{cormain}. 
}
\end{rmk}

%

\section{Characterizations of Purity}\label{sec3}

Let $n, k\in \N$ 
with 
$1\le k\le n$ and let $C$ be an $[n,k]_q$-code. 
We have explained in the Introduction how one can associate an abstract simplicial complex $\Delta_C$ to $C$. 
Note that this complex is independent of the choice of a parity check matrix of $C$. Let $R:= \Fq[x_1, \dots , x_n]$ denote the polynomial ring in $n$ variables over $\Fq$ and let $I_C$ denote the ideal of $R$ generated by the monomials 
$x_{i_1} \cdots x_{i_t}$ where $\{i_1, \dots , i_t\}$ vary over non-faces, i.e., over subsets of $[n]:=\{1, \dots , n\}$ that are not in $\Delta_C$. 
The 
Stanley-Reisner ring $R_C$ corresponding to $\Delta_C$ (with the base field $\Fq$) is, by definition, the quotient $R/I_C$. We call $R_C$ the \emph{Stanley-Reisner ring} associated to $C$. Clearly, $R_C$ is a standard graded $\Fq$-algebra and as noted in 
\cite[\S 1]{GS}, $R_C$ is Cohen-Macaulay and it admits an 
$\N$-graded minimal free resolution 
of the form 
\begin{equation}
\label{eq:resol}
F_k \longrightarrow F_{k-1}\longrightarrow \cdots \longrightarrow F_1\longrightarrow F_0 \longrightarrow R_{\Delta} \longrightarrow 0
\end{equation}
where $F_0= R$ and each $F_i$ is a graded free $R$-module of the form 
\begin{equation}
\label{eq:FiBetaij}
F_i = \bigoplus_{j \in \mathbb{Z}}R(-j)^{\beta_{i,j}} \quad \text{for } i=0, 1, \dots , k.
\end{equation}
The nonnegative integers $\beta_{i,j}$ thus obtained 
are called the \emph{Betti numbers} of~$C$. 
The resolution \eqref{eq:resol} is said to be \emph{pure} of type $(d_0, d_1, \ldots, d_k)$ if for each $i=0,1, \dots , k$, the Betti number $\beta_{i,j}$ is nonzero if and only if $j=d_i$. If, in addition, 
$d_1, \dots , d_k$ are consecutive, 
then the resolution is said to be \emph{linear}. We remark that the Betti numbers 
 $\beta_{i,j}$ as well as the properties of purity and linearity depend only on $C$  and they are independent of the choice of a minimal free resolution of $R_C$. 
 
 The 
 result below is due to Johnsen and Verdure \cite{JV1}; see also \cite[Corollary 3.9]{GS}. 
 
 \begin{prop}\label{MDSlinear}
 Let $C$ be an $[n,k]_q$-code. Then $C$ is an MDS code if and only if $C$ is nondegenerate and every $\N$-graded minimal free resolution of $R_C$ is linear. 
 \end{prop}
 
 We will now recall the intrinsic characterization of purity given in \cite{GS} and alluded to in the Introduction. 
But first, we  review some relevant terminology about codes. 
 
 Let $n,k$ and $C$ be as above. By a \emph{subcode} of $C$ we mean a $\Fq$-linear subspace of $C$. 
 Given a subcode $D$ of $C$, the \emph{support} of $D$ and  the \emph{weight} of $D$ are 
 defined by 
 $$
 \supp(D):=\{ i\in [n] :  \exists \text{  $(c_1, \dots , c_n) \in D$ with $c_i\ne0$}\} \quad \text{and}\quad
 \wt(D):=|\supp(D)|.
 $$
Given any $c\in C$, we often denote  by $\supp(c)$ and $\wt(c)$ the support of $\langle c \rangle$ and the weight of $\langle c \rangle$, respectively, where $\langle c \rangle$ denotes the subcode of $C$ spanned by $c$. 
For $1\le i \le k$, the $i^{\rm th}$ \emph{generalized Hamming weight} of $C$ is defined~by 
$$
d_i(C) := \min\{\wt(D) : D \text{ a subcode of $C$ with } \dim\,D = i\}.
$$
It is well-known that $d_1(C)$ is the minimum distance of $C$ and $d_i(C)< d_{i+1}(C)$ for $1\le i <k$. Note that $C$ is nondegenerate 
if and only if $d_k(C)=n$. 
An $i$-dimensional subcode $D$ of $C$ is said to be \emph{$i$-minimal} if its support is minimal among the supports of all $i$-dimensional subcodes of $C$, i.e., $\supp(D^\prime)\nsubseteq\supp(D)$ for any $i$-dimensional subcode $D'$ of $C$, with $D^\prime\neq D$.

 We are now ready to state (an equivalent version of) the intrinsic characterization of purity given in \cite[Theorem 3.6]{GS}. 

\begin{prop}\label{1}
Let $C$ be an $[n,k]_q$-code and let $d_1<\cdots < d_k$ be its generalized Hamming weights. Also, let $R_C$ be the Stanley-Reisner ring associated 
to $C$. Then every $\N$-graded minimal free resolution of $R_C$ is not pure if and only if  there exists an $i\in \{1,\ldots,k\}$ and an $i$-minimal subcode $D_i$ of $C$ such that $\wt(D_i) > d_i$. 
\end{prop}

We summarize below the results in \cite{GS} about the purity and non-purity of graded minimal free resolutions of 
Stanley-Reisner ring associated to Reed-Muller codes. 

\begin{prop}\label{mainGS}
Let $m, r$ be integers such that $m\ge 1$ and $0\le r \le m(q-1)$. Also, let $t, s$ be unique nonnegative integers satisfying \eqref{eq:ts}. Then every $\N$-graded minimal free resolution of the
Stanley-Reisner ring associated to ${\mathcal RM}_q(r,m)$ is 
\begin{enumerate}
\item[{\rm (i)}] pure if $r=1$, 
\item[{\rm (ii)}] not pure if $q=2$, $m \ge 4$, and $1< r \le m-2$, and 
\item[{\rm (iii)}] not pure if $m \ge 2$, $1<r<m(q-1)-1$, and $s\ne 1$. 
\end{enumerate}
\end{prop}

\begin{proof}
The assertion in (i) is proved in \cite[Theorem 4.1]{GS}, while the assertions in (ii) and (iii) are proved in \cite[Proposition 4.4]{GS} and \cite[Theorem 4.11]{GS}, respectively. 
\end{proof}

The values of $q, m,r$ 
not covered by (i)--(iv) in Lemma~\ref{lem2} and (i)--(iii) in Proposition~\ref{mainGS} are precisely $q \ge 3$, $m\ge 2$, and $r=q, \; 2q-1,  \dots , (m-1)q -  (m-2)$, except that $(m-1)q -  (m-2)$ is excluded if $q=3$. This is taken care of by the following. 

\begin{lemma}\label{main}
 
Let $m,r$ be 
integers 
such that $m \ge 2$ and $1<r<m(q-1)-1$. Also let $t,s$ be unique integers satisfying \eqref{eq:ts}. 
Assume that $q\ge 3$ and also that  
$s=1$. 
Then 
every $\N$-graded minimal free resolution of the
Stanley-Reisner ring associated to the Reed-Muller code ${\mathcal RM}_q(r,m)$ is not pure. 
\end{lemma}

\begin{proof}
The conditions on $m, r$ and our assumptions imply that $1\le t \le m-1$ and moreover if $q=3$, then $1\le t \le m-2$. 
Also note that by Proposition~\ref{RM-min}, the minimum distance of ${\mathcal RM}_q(r,m)$ is given by $d= (q-1)q^{m-t-1}$.
We will divide the proof in two cases according as $q> 3$ and $q=3$. 

\smallskip

{\bf Case 1.} $q>3$. 

Write $\F_q=\{\omega_1,\ldots,\omega_q\}$, and let $\omega'_1, \omega'_2$ be two distinct elements of $\F_q$. Define 
\begin{eqnarray*}
Q(X_1,\ldots,X_m) := \left(\prod_{i=1}^{t-1}(X_i^{q-1}-1)\right)\left(\prod_{j=3}^{q}(X_t-\omega_j)\right)\left(\prod_{k=1}^{2}(X_{t+1}-\omega'_k)\right).
\end{eqnarray*}
Then ${\rm deg}(Q)=(t-1)(q-1)+(q-2)+2=(t-1)(q-1)+q=t(q-1)+1=r$, and thus $Q\in V_q(r,m)$. For $i=1,2$, let 
\begin{eqnarray*}
A_i := \left\{{\bf a} =(a_1, \dots , a_m) \in \F_q^m:a_1=\cdots = a_{t-1}=0,\; a_t=\omega_i \text{ and }  a_{t+1}\not\in \{\omega'_1,\omega'_2\} \right\}. 
\end{eqnarray*}
Then ${\rm Supp}(c^{ }_Q)=A_1\cup A_2$. Observe that $A_1$ and $A_2$ are disjoint. Consequently, 
$$
{\rm wt}(c^{ }_Q)=2(q-2)q^{m-t-1} \quad \text{and therefore} \quad \wt(c^{ }_Q) >d  = (q-1)q^{m-t-1},
$$ 
where the last inequality follows since $q>3$. 
Thus $c^{ }_Q$ is not a minimum weight codeword. 
We claim that the $1$-dimensional subcode $\langle c^{ }_Q \rangle$ is $1$-minimal. This claim together with  Proposition~\ref{1} would imply the desired result. To prove the claim, assume the contrary. Thus, suppose there is $F\in V_q(r,m)$, such that $c^{ }_F$ is a minimum weight codeword of $\mathcal{RM}_q(r,m)$ and ${\rm Supp}(c^{ }_F)\subsetneq {\rm Supp}(c^{ }_Q)$. 
By Proposition~\ref{RM-min}, 
$F$ must be of the form
\begin{equation}\label{eq:minF}
F(X_1,\ldots,X_m)= \omega_0\left(\prod_{i=1}^{t}(1-L_i^{q-1})\right)(L_{t+1}-\omega)
\end{equation}
for some linearly independent linear polynomials $L_1,\ldots,L_{t+1}$ in $\Fq[X_1, \dots , X_m]$ and some $\omega_0,\omega\in \F_q$ with 
$\omega_0\neq 0$. Note that ${\rm Supp}(c^{ }_F) = A'$, where 
\begin{equation}\label{eq:Aprime}
A':= \left\{ {\bf a}=(a_1,\ldots,a_m)\in \F_q^m : L_i({\bf a})=0 \text{ for $1\le i \le t$ and }
L_{t+1}({\bf a})\neq \omega\right\}. %
\end{equation}
Since ${\rm Supp}(c^{ }_F)\subset {\rm Supp}(c^{ }_Q) $, we obtain 
$A'\subset A_1\cup A_2$. 
We now assert that $A'$ is disjoint from one of the $A_i$. Indeed, if the assertion were not true, then we can choose $P_i\in A'\cap A_i$ for $i=1,2$. Write $b_i := L_{t+1}(P_i)$ for $i=1,2$.
Since $P_i\in A'$, we see that $b_i\neq \omega$ for $i=1,2$. Now pick $\lambda\in \F_q$ such that $\lambda\neq 0,1$  and $(1-\lambda)b_1+\lambda b_2\neq \omega$, which is possible because $q\geq 4$.\footnote[1]{If $b_1=b_2$, then the only condition on $\lambda$ is that $\lambda\neq 0,1$, whereas if $b_1\ne b_2$, then it suffices to choose $\lambda \in \Fq$ such that $\lambda\neq 0,1$ and $\lambda\neq (\omega - b_1)/(b_2-b_1)$.}
Define $P_\lambda:=(1-\lambda)P_1+\lambda P_2$. 
Then $P_\lambda\in A'$, and this contradicts the inclusion $A'\subset A_1\cup A_2$ because the $t^{\rm th}$ coordinate of $P_\lambda$ is neither $\omega_1$ nor $\omega_2$. This proves the above assertion. Thus ${\rm Supp}(c^{ }_F)=A'\subseteq A_i$ for some $i$. But then 
$(q-1)q^{m-t-1}\leq (q-2)q^{m-t-1}$, which is a contradiction. This proves the claim and hence the desired result 
when $q>3$.
\smallskip

{\bf Case 2.} $q=3$.

In this case $1\le t \le m-2$,  
as noted earlier.  
Write $\Fq = \{\omega_1,\omega_2,\omega_3\}$. 
Define
$$
Q(X_1,\ldots,X_m) := \bigg(\prod_{i=1}^{t-1}(X_i^{q-1}-1)\bigg)(X_t-\omega_3)(X_{t+1}-\omega_3)(X_{t+2}-\omega_3). 
$$
Then ${\rm deg}(Q)=(t-1)(q-1)+3=t(q-1)+1=r$, since $q=3$, and so $Q\in V_q(r,m)$. 
Let $E:= \left\{{\bf a} =(a_1, \dots , a_m) \in \F_q^m: a_1=\cdots=a_{t-1}=0\right\}$, and  for $i=1,2$, 
let 
\begin{eqnarray*}
A_i &:=& \left\{{\bf a} =(a_1, \dots , a_m) \in E :  a_t=\omega_i \text{ and } a_{t+1}, a_{t+2} \in \{\omega_1,\omega_2\} \right\}, \\
A'_i &:=& \left\{{\bf a} =(a_1, \dots , a_m) \in E :  a_{t+1}=\omega_i \text{ and } a_{t}, a_{t+2} \in \{\omega_1,\omega_2\} \right\},  \text{ and}\\
A''_i &:=& \left\{{\bf a} =(a_1, \dots , a_m) \in E :  a_{t+2}=\omega_i \text{ and } a_{t}, a_{t+1} \in \{\omega_1,\omega_2\} \right\}. 
\end{eqnarray*}
Then ${\rm Supp}(c^{ }_Q)=A_1\cup A_2=A'_1\cup A'_2=A''_1\cup A''_2$ and $\wt(c^{ }_Q)= 2^3 q^{m-t-2}$. 
Note that $\wt(c^{ }_Q) >(q-1)q^{m-t-1}$, since $q=3$. 
Thus, as in Case 1, it suffices to show  that there does not exist any $F\in V_q(r,m)$ such that $c^{ }_F$ is a minimum weight codeword and ${\rm Supp}(c^{ }_F)\subsetneq {\rm Supp}(c^{ }_Q)$. Suppose, if possible, there is such $F$. Then it must be of the form \eqref{eq:minF}, and its support is given by the set $A'$ in \eqref{eq:Aprime}. Now write 
$\Fq\setminus \{\omega\} = \{u_1, u_2\}$, and for $i=1,2$, let 
$$
B_i := \left\{{\bf a}=(a_1,\ldots,a_m)\in \F_q^m : L_i({\bf a})=0 \text{ for } 1\le i \le t \text{ and } L_{t+1}({\bf a})=u_i\right\}.
$$
Note that each $B_i$ is an affine space (i.e., a translate of a linear subspace) in $\Fq^m$ and $\supp(c^{ }_F) = B_1\cup B_2$.  Thus $B_1 \cup B_2 \subset A_1 \cup A_2$. We claim that $B_1 \subseteq A_i$ for some $i\in \{1, 2\}$. Indeed, if this were not true, then we can find $P_i \in B_1\cap A_i$ for each $i=1, 2$. Since $q=3$, we can choose $\lambda \in \F_q$ such that $\lambda\ne 0, 1$. Consider  $P_\lambda := (1-\lambda)P_1+\lambda P_2$. 
Since $B_1$ is an affine space, $P_\lambda \in B_1$. On the other hand, the $t^{\rm th}$ coordinate of $P_\lambda$ is neither $\omega_1$ nor $\omega_2$, and hence $P_\lambda\not\in A_1\cup A_2$. This  contradicts the inclusion $B_1  \subset A_1 \cup A_2$, and so the Claim is proved. In a similar manner, we see that $B_1 \subseteq A'_j$ 
and $B_1 \subseteq A''_k$ 
for some $j, k\in \{1, 2\}$. It follows that $B_1 \subseteq A_i \cap A'_j \cap A''_k$. But clearly, $|B_1| = q^{m-t-1}$ 
and $ |A_i \cap A'_j \cap A''_k|= q^{m-t-2}$. So we obtain 
$q^{m-t-1}\leq q^{m-t-2}$, which is a contradiction. This completes the proof.
\end{proof}

We are now ready to prove the main result of this article. 

\begin{thm}\label{mainthm}
Let $m, r \in \N$ be 
such that $m\ge 1$ and $0\le r \le m(q-1)$. Then every $\N$-graded minimal free resolution of the Stanley-Reisner ring associated to the Reed-Muller code $\mathcal{RM}_q(r,m)$ is pure if and only if 
$m=1$ or $r\le 1$ or $r \ge m(q-1) -1$. 
\end{thm}

\begin{proof}
Follows from Lemma~\ref{lem2}, Proposition~\ref{MDSlinear}, Proposition~\ref{mainGS}, and Lemma~\ref{main}.
\end{proof}

As an application, we show that the converse of the result in Lemma~\ref{lem2} is true. 

\begin{cor}\label{cormain}
Let $m, r \in \N$ be 
such that $m\ge 1$ and $0\le r \le m(q-1)$. Then the Reed-Muller code $\mathcal{RM}_q(r,m)$ is an MDS code if and only if 
$m=1$ or $ r = 0$ or $r \ge m(q-1)-1$. 
\end{cor}

\begin{proof}
If $m=1$ or $ r = 0$ or $r \ge m(q-1)-1$, then by  Lemma~\ref{lem2}, $\mathcal{RM}_q(r,m)$ is an MDS code. 
Conversely, suppose  $\mathcal{RM}_q(r,m)$ is an MDS code. Then by  Proposition~\ref{MDSlinear}, every $\N$-graded minimal free resolution of its Stanley-Reisner ring is pure. So by Theorem~\ref{mainthm}, we must have $m=1$ or $ r \le 1$ or $r \ge m(q-1)-1$. If $m\ge 2$, then the case $r=1$ is ruled out because by \cite[Theorem~4.1]{GS}, the generalized Hamming weights (which coincide with the ``shifts'' in the resolution) of  $\mathcal{RM}_q(1,m)$ are given by $d_i = q^m-\lfloor q^{m-i}\rfloor$ for $1\le i\le m+1$, and these are clearly non-consecutive if $m\ge 2$, and so by Proposition~\ref{MDSlinear}, $\mathcal{RM}_q(1,m)$ cannot be an MDS code  if $m\ge 2$. Thus we must have  
$m=1$ or $ r = 0$ or $r \ge m(q-1)-1$. 
\end{proof}


\end{document}